\documentclass[12pt]{article}
\usepackage{amsmath,amssymb,amsthm}
\usepackage{color}
\numberwithin{equation}{section}
\newtheorem{thm}{Theorem}[section]

\def\address#1#2{\begingroup
\noindent\parbox[t]{16cm}{%
\small{\scshape\ignorespaces#1}\par\vskip1ex
\noindent\small{\itshape E-mail address}%
\/: #2\par\vskip4ex}\hfill%
\endgroup}%

\title{Escape rate of the Brownian motions on hyperbolic spaces}
\author{Yuichi Shiozawa\footnote{Supported in part by the Grant-in-Aid for Scientific Research (C) 26400135.}}
\begin{document}
\maketitle 

\begin{abstract}
We discuss the escape rate of the Brownian motion on a hyperbolic space. 
We point out that the escape rate is determined by using the Brownian expression of the radial part  
and a generalized Kolmogorov's test for the one dimensional Brownian motion.
\end{abstract}

\section{Introduction}
Let ${\mathbb H}^d$ be the $d$-dimensional hyperbolic space
and ${\mathbf M}=(\{X_t\}_{t\geq 0}, \{P_x\}_{x\in {\mathbb H}^d})$ the Brownian motion 
on ${\mathbb H}^d$ generated by the half of the Laplace-Beltrami operator.
For a fixed point $o\in {\mathbb H}^d$, define $P=P_o$ and $R_t=d(o,X_t)$. 
In this note, we show

\begin{thm}\label{criteria}
Let $g(t)$ be a positive function on $(0,\infty)$ such that for some $t_0>0$, 
$\sqrt{t}g(t)$ is nondecreasing and $g(t)/\sqrt{t}$ is bounded for all $t\geq t_0$.  
\begin{enumerate}
\item For the function $r_1(t):=(d-1)t/2+\sqrt{t}g(t)$,  
\begin{equation}\label{upper}
P\left(\text{there exists $T>0$ such that 
$R_t<r_1(t)$ for all $t\geq T$}\right)=1 \ \text{or $0$}
\end{equation}
according as 
\begin{equation}\label{kol-test}
\int_{\cdot}^{\infty}\left(1\vee g(t)\right)e^{-g(t)^2/2}\frac{{\rm d}t}{t}<\infty \ \text{or $=\infty$}.
\end{equation}
\item For the function $r_2(t):=(d-1)t/2-\sqrt{t}g(t)$, 
\begin{equation}\label{lower}
P\left(\text{there exists $T>0$ such that $R_t>r_2(t)$ for all $t\geq T$}\right)=1 \ \text{or $0$}
\end{equation}
according as {\rm (\ref{kol-test})} holds.
\end{enumerate}
\end{thm}

The function $r_1(t)$ is called an {\it upper rate function} for ${\mathbf M}$ 
if the probability in (\ref{upper}) is $1$. 
By the same way, the function $r_2(t)$ is called a {\it lower rate function} for ${\mathbf M}$ 
if the probability in (\ref{lower}) is $1$. 
According to Theorem \ref{criteria}, we have for $c>0$,
\begin{itemize}
\item the function $r(t):=(d-1)t/2+\sqrt{ct\log\log t}$ 
is an upper rate function for ${\mathbf M}$ if and only if $c>\sqrt{2}$;
\item the function $r(t):=(d-1)t/2-\sqrt{ct\log\log t}$ is a lower rate function for ${\mathbf M}$
if and only if $c>\sqrt{2}$.
\end{itemize}

For the Brownian motions on Riemannian manifolds, more 
generally symmetric diffusion processes, 
upper and lower rate functions are given in terms of the volume growth rate 
(\cite{BS05, G99-0, G99, GH09, HQ10, O16}). 
As for the upper rate functions, the results in 
\cite{G99-0, G99, GH09, HQ10, O16} 
are applicable to the Brownian motions 
on Riemannian manifolds with exponential volume growth rate, as so to ${\bf M}$; 
however, as for the lower rate functions, the results in \cite{BS05, G99-0, G99} 
are not applicable to ${\bf M}$ 
because the doubling condition is imposed on the volume growth.
Grigor'yan and Hsu \cite{GH09} also discussed the sharpness of the upper rate functions 
for ${\mathbf M}$ or 
for the Brownian motion on  a model manifold, that is, a spherically symmetric Riemannian manifold with a pole.
Using the fact that 
\begin{equation}\label{rad-lim}
\lim_{t\rightarrow\infty}\frac{R_t}{t}= \frac{d-1}{2}, \quad \text{$P$-a.s.}
\end{equation}
(which follows from (\ref{rad-exp}) below), 
they remarked that 
the function $r(t)=ct$ is an upper rate function for ${\bf M}$ if $c>(d-1)/2$, 
and not if $0<c<(d-1)/2$. 
This observation is still valid for the lower rate functions. 
See also \cite{I88} for the result of the law of the iterated logarithms-type 
to the Brownian motions on model manifolds.

For the proof of Theorem \ref{criteria}, 
we make use of the Brownian expression of the radial part $R_t$ ((\ref{rad-exp}) below) as in 
\cite{GH09, I88}, 
together with a generalized version of Kolmogorov's test 
for the one dimensional Brownian motion 
(\cite{K97, K98}). 
In fact, the integral in (\ref{kol-test}) is the same with that in this test.
The assumption on $g(t)/\sqrt{t}$ will be needed  
in (\ref{comp-int-test}) and (\ref{comp-int-test-2}) below.

\section{Proof of Theorem \ref{criteria}}

Let ${\mathbf B}=(\{B_t\}_{t\geq 0},P)$ be the one dimensional Brownian motion starting from the origin. 
Then a generalized Kolmogorov's test holds: 
\begin{thm}\label{thm-g-kol}{\rm (\cite[Theorem 3.1 and Lemma 3.3]{K97} and \cite[Theorem 2.1]{K98})}
Under the full conditions of Theorem {\rm \ref{criteria}},  
$$P\left(\text{there exists $T>0$ such that $|B_t|<\sqrt{t}g(t)$ for all $t\geq T$}\right)=1 \ \text{or $0$}$$
according as {\rm (\ref{kol-test})} holds. 
This assertion is valid even if $|B_t|$ in the equality above is replaced by $B_t$ or $-B_t$.
\end{thm}

By comparison with Kolmogorov's test (see, e.g., \cite[4.12]{IM74}), 
we do not need to assume that $g(t)\nearrow\infty$ as $t\rightarrow\infty$ in Theorem \ref{thm-g-kol}.
\medskip

{\it Proof of Theorem {\rm \ref{criteria}}.} 
Recall that  ${\mathbf M}=(\{X_t\}_{t\geq 0}, P)$ is the Brownian motion 
on ${\mathbb H}^d$ starting from a fixed point $o\in {\mathbb H}^d$ 
and $R_t=d(o,X_t)$ is the radial part of $X_t$.  
Then by \cite[Example 3.3.3]{H02},
\begin{equation}\label{rad-exp}
R_t=B_t+\frac{d-1}{2}\int_0^t\coth R_s\,{\rm d}s.
\end{equation}

Assume the full conditions of Theorem \ref{criteria}. 
We first discuss the lower bound of $R_t$. 
Since $\coth x\geq 1$ for any $x>0$, we obtain by (\ref{rad-exp}),
$$R_t\geq B_t+\frac{d-1}{2}t \quad \text{for any $t\geq 0$}.$$
Hence if the integral in (\ref{kol-test}) is convergent, 
then the probability in (\ref{lower}) is $1$ by Kolmogorov's test. 
By the same way, if the integral in (\ref{kol-test}) is divergent, 
then the probability in (\ref{upper}) is $0$. 

We next discuss the upper bound of $R_t$. 
Since ${\bf M}$ is transient (e.g., see \cite[Subsection 3.2]{G99}), 
we have $P(\lim_{t\rightarrow\infty}R_t=\infty)=1$ so that $P(\lim_{t\rightarrow\infty}\coth R_t=1)=1$. 
Therefore,
$$\lim_{t\rightarrow\infty}\frac{1}{t}\int_0^t\coth R_s\,{\rm d}s=1, \quad \text{$P$-a.s.}$$
By noting that $B_t=o(t)$ as $t\rightarrow\infty$, there exists $c>0$ such that 
$P(A)=1$ for 
$$A:=\left\{\text{there exists $T_1>0$ such that $R_t \geq c t$ for all $t\geq T_1$}\right\}.$$ 

Under the event $A$,  
$$\coth R_s-1=\frac{2}{e^{2R_s}-1}\leq \frac{2}{e^{2cs}-1} \quad \text{for any $s\geq T_1$},$$
which implies that for all $t\geq T_1$,
\begin{equation*}
\begin{split}
\int_0^t(\coth R_s-1)\,{\rm d}s
&=\int_0^{T_1}(\coth R_s-1)\,{\rm d}s+\int_{T_1}^t(\coth R_s-1)\,{\rm d}s\\
&\leq \int_0^{T_1}(\coth R_s-1)\,{\rm d}s+\int_{T_1}^{\infty}\frac{2}{e^{2cs}-1}\,{\rm d}s
=:C_{T_1}.
\end{split}
\end{equation*}
Since there exists  an integer valued random variable $N$ such that 
\begin{equation}\label{n}
\frac{d-1}{2}C_{T_1}\leq N,
\end{equation}
we obtain for such $N$, 
\begin{equation}\label{upper-bound}
\begin{split}
R_t
&=B_t+\frac{d-1}{2}t+\frac{d-1}{2}\int_0^t(\coth R_s-1)\,{\rm d}s\\
&\leq B_t+\frac{d-1}{2}t+\frac{d-1}{2}C_{T_1}
\leq B_t+\frac{d-1}{2}t+N \quad \text{for all $t\geq T_1$}.
\end{split}
\end{equation}

Assume first that the integral in (\ref{kol-test}) is convergent. 
Then 
there exists a positive constant $c_n$ for each  $n\geq 1$ such that 
the function $h_1^{(n)}(t):=g(t)-n/\sqrt{t}$ satisfies 
\begin{equation}\label{comp-int-test}
\begin{split}
\int_{\cdot}^{\infty} \left(1\vee h_1^{(n)}(t)\right)
e^{-h_1^{(n)}(t)^2/2}\frac{{\rm d}t}{t}
&\leq c_n\int_{\cdot}^{\infty} \left(1\vee g(t)\right) e^{-g(t)^2/2}\frac{{\rm d}t}{t}<\infty.
\end{split}
\end{equation}
Hence Theorem \ref{thm-g-kol} implies that for each $n\geq 1$, 
$$P\left(\text{there exists $T>0$ such that $|B_t|<r_1^{(n)}(t)$ for all $t\geq T$}\right)=1$$
for $r_1^{(n)}(t):=\sqrt{t}h_1^{(n)}(t) (=\sqrt{t}g(t)-n)$. 
In particular, we get $P(B_1)=1$ for 
$$B_1:=\left\{\text{for each $n\geq 1$, there exists $S_n>0$ such that $|B_t|<r_1^{(n)}(t)$ for all $t\geq S_n$}\right\}.$$
Under the event $A\cap B_1$, 
since 
there exists $T_2>0$ for $N\geq 1$ in (\ref{n}) such that 
$$B_t<r_{1}^{(N)}(t)=\sqrt{t}g(t)-N \quad \text{for all $t\geq T_2$},$$
we have  by (\ref{upper-bound}),
$$R_t<\frac{d-1}{2}t+\sqrt{t}g(t) \quad \text{for all $t\geq T_1\vee T_2$}.$$
Therefore, the probability in (\ref{upper}) is $1$.

Assume next that the integral in (\ref{kol-test}) is divergent. 
Then by the same way as in  (\ref{comp-int-test}), 
the function $h_2^{(n)}(t):=g(t)+n/\sqrt{t}$ satisfies for each $n\geq 1$,
\begin{equation}\label{comp-int-test-2}
\int_{\cdot}^{\infty} \left(1\vee h_2^{(n)}(t)\right)
e^{-h_2^{(n)}(t)^2/2}\frac{{\rm d}t}{t}=\infty.
\end{equation}
Hence Theorem \ref{thm-g-kol} yields that  for each $n\geq 1$, 
$$P\left(\text{for any $t>0$, there exists $T\geq t$ such that $B_T\leq -r_2^{(n)}(T)$}\right)=1$$
for $r_2^{(n)}(t):=\sqrt{t}h_2^{(n)}(t) (=\sqrt{t}g(t)+n)$. 
In particular, $P(B_2)=1$ for 
$$B_2:=\left\{\text{for each $n\geq 1$, there exists $U_n\geq t$ for any $t>0$ 
such that $B_{U_n}\leq -r_2^{(n)}(U_n)$}\right\}.$$
Under the event $A\cap B_2$, 
since there exists $T_3\geq t\vee T_1$ for any $t>0$ and $N\geq 1$ in (\ref{n}) such that 
$$B_{T_3}\leq -r_2^{(N)}(T_3)=-\sqrt{T_3}g(T_3)-N,$$
we have for such $T_3$,
$$R_{T_3}\leq \frac{d-1}{2}T_3-\sqrt{T_3}g(T_3)$$
by (\ref{upper-bound}).  
Therefore, the probability in (\ref{lower}) is $0$.
\qed
\bigskip

\noindent
{\bf Acknowledgements.} \ 
The author would like to thank Professor Masayoshi Takeda 
for his valuable discussion motivating this work 
and his comment on the draft of this paper.

\address{
Yuichi Shiozawa\\
Graduate School of Natural Science and Technology \\
Department of Environmental and Mathematical Sciences \\
Okayama University \\
Okayama 700-8530,
Japan
}
{shiozawa@ems.okayama-u.ac.jp}
\end{document}